\documentclass[a4paper,10pt]{article}

\usepackage{amsmath,amssymb,amsthm,amsopn,latexsym,amstext}
\usepackage{color}
\usepackage{hyperref}
\usepackage[margin=1in]{geometry}

\usepackage{comment}

\usepackage{numbysec}

\usepackage[normalem]{ulem}

\usepackage[T1]{fontenc}
\usepackage[cp1250]{inputenc}

\newtheorem{thm}{Theorem}[section] 
\newtheorem{lem}[thm]{Lemma}
 \newtheorem{prop}[thm]{Proposition}
\newtheorem{cor}[thm]{Corollary} 

\newtheorem{example}{Example}
\newtheorem{remark}{Remark}

\DeclareMathOperator{\diver}{div}

\def \EE {\mathbb{E}} \def \PP {\mathbb{P}} \def \RR {\mathbb{R}} \def
\NN {\mathbb{N}} \def \Rd {{\RR^d}}       
\newcommand{\set}[1]{\left\{#1\right\}} \def \eps {\varepsilon}

        \def \CIII
{C}

\def \tp{\tilde{p}}
\def \tEE {\tilde{\mathbb{E}}} 
\def \tPP {\tilde{\mathbb{P}}}

\def \ga{\gamma}

\def \tG {\tilde G}

\def \tP {\tilde P}

\def \de {\delta}

\def \al {\alpha}
\def \la{\lambda}

\def \jo {\iota}
\def \ka {\kappa}

\def \Lt {L^2(D) }
\def\bbR{\mathbb R}

\def\bbE{\mathbb E}
\def\proof{\noindent{\em Proof. }}

\def\u{\underline}

\def \ka{\kappa}

\definecolor{tk}{rgb}{0.4,0,0.6}

\definecolor{kb}{rgb}{0.2,0.7,0.1}

\title{Principal eigenvalue of the fractional Laplacian with a large incompressible drift
\footnote{
The research of K.B. was partially supported by NCN grant 2012/07/B/ST1/03356.
T.K. acknowledges
the support of the Polish Ministry of Higher Education
grant NN201419139. }
}
\author{Krzysztof Bogdan \and Tomasz Komorowski \footnote{Institute of Mathematics of the
Polish Academy of Sciences, ul. \'Sniadeckich 8, 00-956 Warszawa, Poland, e-mails: bogdan@pwr.wroc.pl, komorow@hektor.umcs.lublin.pl}}

\date{\today}

\numberbysection

\begin{document}
\maketitle

\begin{abstract} 
We study the  principal Dirichlet eigenvalue of the operator $L_A=\Delta^{\al/2}+Ab(x)\cdot\nabla$,  on a bounded $C^{1,1}$ regular domain $D$.
Here
$\al\in(1,2)$, 
$\Delta^{\alpha/2}$ is the fractional 
Laplacian, $A\in\bbR$, and $b$
is a 
bounded $d$-dimensional 
divergence-free
vector field
in the Sobolev space 
$W^{1,2d/(d+\alpha)}(D)$.
We prove that the  eigenvalue remains bounded, as $A\to+\infty$, if and only if 
$b$
has
non-trivial first integrals
in the domain of the quadratic
form 
of 
$\Delta^{\alpha/2}$ 
for the Dirichlet condition.
\end{abstract}

\section{Introduction}

This article is motivated by 
the following
result 
of  Berestycki, et al.  given in
\cite{MR2140256}
for the Laplacian
perturbed by a divergence-free drift  in dimensions  $d\ge 2$.
Let $D\subset\Rd$ be a bounded $C^2$ regular open set and  let  
$b(x)=(b_1(x),\ldots,b_d(x)):\Rd\to\Rd$ be a bounded $d$-dimensional vector field  
such that 
$
\diver b=0$ on $D$
in the sense of distributions (distr.), i.e. 
\begin{equation}
\label{div-free}
\int_{\Rd}b(x)\cdot\nabla\phi(x)dx=0,\quad 
\phi\in C_c^\infty(D). 
\end{equation}
For $A\in\bbR$,  let ($\phi_A,\lambda_A$) be the principal eigen-pair corresponding to the Dirichlet problem for the operator
$\Delta + Ab(x)\cdot\nabla$. 
 Theorem 0.3 of \cite{MR2140256} asserts that $\lambda_A$ remains bounded as $A\to+\infty$, if and only if the  equation 
\begin{equation}
\label{int-first}
{\rm div}\left(wb\right) 
=0\quad \mbox{ (distr.) on $D$}
\end{equation}
has 
 a solution $w$ (called a first integral of $b$), such that $w\neq 0$ and
$w\in H_0^1(D)$.

The result 
can be interpreted intuitively
in the following way: 
functions $w$ satisfying \eqref{int-first}
are constant along the 
flow of the vector field $Ab(x)$ (see Section~\ref{sec:ub}), and
the existence of (non-trivial) first integrals allows for flow lines that are contained in $D$.
On the other hand, if no such $w$ exist, then the flow leaves $D$ with 
speed proportional to $A$.
Adding the Laplacian $\Delta$ to $b\cdot \nabla$, or equivalently the Brownian motion 
to the flow,
results in a stochastic process 
whose trajectories 
gradually depart from the integral curves of $b$,
but the general picture is similar:
if 
nontrivial first integrals 
exist,
then the
trajectories may 
remain
in 
$D$ 
with positive probability 
during a 
finite time interval, even as $A\to+\infty$. 
In 
this case we are lead to a nontrivial 
limiting transition mechanism 
between 
the flow lines.
The result described in the foregoing 
enjoys many extensions and has proved  
quite useful in various 
applications describing the influence of a fluid flow on a diffusion,  see for example 
\cite{MR2753290,MR2434887,MR2677980,MR2748635}. In the context of
  a compact, connected Riemannian manifold a sufficient and necessary condition for $\lambda_A$ to
  remain bounded, as $A\to+\infty$, expressed  in terms of the eigenspaces of the
advection  operator $b(x)\cdot \nabla$, has been given in \cite[Theorem 1]{MR2563731}.

The purpose of the present paper is to  
verify
that a
similar property of the principal eigenvalue holds
when the classical Laplacian is replaced by the fractional Laplacian $\Delta^{\alpha/2}$
with $\alpha\in(1,2)$. We consider ${\cal I}_0^{\al}$ defined as the set of all the
nonzero first integrals
in the Sobolev space $H^{\al/2}_0(D)$ equipped with the norm
coming from 
the 
Dirichlet form ${\cal E}^\alpha$
of
$\Delta^{\alpha/2}$
(see
\eqref{011109b} below). 
The Sobolev norm condition on the first integrals reflects smoothing properties of the Green function of the fractional Laplacian, while \eqref{int-first} is related to the flow defined by $b$.

The 
main difficulty 
in 
our development
stems from roughness of general elements of $H^{\alpha/2}_0(D)$
and  non-locality of
$\Delta^{\alpha/2}$,
which prevent us  from a direct application of the 
differential calculus
in the way it has been done in \cite{MR2140256}.
Instead, we use conditioning suggested by a paper of Bogdan and Dyda \cite{MR2663757}, approximation techniques for flows given by DiPerna and Lions in \cite{MR1022305},
and the 
properties of the Green function 
and 
heat kernel
of
gradient perturbations of
$\Delta^{\alpha/2}$ obtained by Bogdan, Jakubowski  in \cite{MR2892584} and Chen, et al. in \cite{2010arXiv1011.3273C}  for $\al\in(1,2)$ and bounded $C^{1,1}$-regular 
open sets $D$. 
These 
properties
allow to 
define and study, via  
the classical Krein-Rutman theorem and compactness arguments,
the principal eigen-pair $(\lambda_A,\phi_A)$ for 
$L_A=\Delta^{\al/2}+Ab\cdot\nabla$ and
$\al\in(1,2)$. 
Our main result can be stated as follows.
\begin{thm}
\label{main-thm}
Suppose that $D\subset \Rd$ is 
a bounded domain with 
$C^{1,1}$-regular 
boundary that is simply connected, i.e.  $\bar D^c$ - the complement
of $\bar D$ - is connected. 
If $\alpha\in(1,2)$, and $b\in L^\infty(D)\cap
W^{1,2d/(d+\alpha)}(D)$ is of zero divergence, then 
\begin{equation}
\label{081007}
\lim_{A\to+\infty}\lambda_A={
\inf\limits_{w\in {\cal I}_0^{\al},\|w\|_2=1}}
{\cal E}^\alpha(w,w),
\end{equation}
and the infimum is attained. Here we
use the convention that $\inf \emptyset=+\infty$,
hence
$\lim_{A\to+\infty}\lambda_A=+\infty$ if and only if the zero function is the only first integral.
\end{thm}
Equality \eqref{081007} results from
the following lower and upper bounds
of $\lambda_A$,
\begin{equation}
\label{081007a}
\liminf_{A\to+\infty}\lambda_A\ge \inf_{w\in {\cal I}^{\al},\|w\|_2=1}{\cal E}^{\al}(w,w),
\end{equation}
\begin{equation}
\label{091007b}
\sup_{A\in \RR} \lambda_A\le \inf_{w\in {\cal I}^{\al},\|w\|_2=1}{\cal E}^{\al}(w,w).
\end{equation}
The bounds are proved in  
Sections~\ref{sec:lb} and 
\ref{sec:ub}, correspondingly.
In Section
\ref{proof-coro}
 we explain 
that the
minimum on the right hand side of \eqref{081007} is attained,
and we finish
the proof of the theorem.

Comparing our approach with the arguments used in the case of 
local operators,
  cf. \cite{MR2140256,MR2563731}, we note that the use of the Green function seems
  more robust 
whenever we lack 
sufficient
  differentiability of functions appearing in variational
  formulas. 
Recall that in the present case 
we need to deal with $H^{\al/2}_0(D)$,
which limits the applicability of the
  arguments based on the 
usual 
  differentiation rules of
the classical calculus,
e.g. the Leibnitz formula or the chain rule.  
We consider the use of the Green function as one of the major features of our approach. 
In addition, the non-locality of the
  quadratic forms  forces a substantial modifications of several other arguments,
e.g. those involving
conditioning of
nonlocal operators and quadratic forms in the proof of the
upper bound \eqref{091007b} in Section \ref{sec:ub}.
 Finally, we stress the fact that the Dirichlet fractional Laplacian on a bounded domain $D$ {\it is not} a fractional power of the Dirichlet Laplacian on $D$, e.g. the eigenfunctions of these operators 
have a different power-type decay at the boundary, see \cite{MR2876409,MR2217951,MR2415141} in this connection.

As a preparation for the proof, we recall  in Section \ref{sec:P}  the
estimates of \cite{MR2892584,2010arXiv1011.3273C} for the Green function and transition 
density of $L_A$ for the Dirichlet problem on $D$.
{These functions are defined using Hunt's formula \eqref{eq:Hunt}, which in principle requires the drift $b(x)$ to be defined on the entire $\bbR^d$. We show however,
in Corollary \ref{cor010212}, that they are 
determined by the restriction of the drift to the domain $D$.
In Section \ref{sec3} we prove that the corresponding Green's and transition operators are compact, see 
{Lemmas} \ref{lem:GDc1} and \ref{lem:GDc1KB}.   This result is used to define the principal eigen-pair of $L_A$, via the Krein-Rutman theorem.}
In Theorem \ref{thm011307} of Section \ref{sec4}  
we prove that the domains of $\Delta^{\alpha/2}$ and  $L_A$ in $L^2(D)$
coincide. 
In Section~\ref{sec:pt11} we employ the 
bilinear 
form of $L_A$ to estimate the principal eigenvalue. The technical assumption
$\nabla b\in L^{2d/(d+\alpha)}(D)$
is only 
needed in Sections~\ref{sec:ub}
to characterize the first integrals of $b$ by means of
the theory of flows developed by DiPerna and Lions
in \cite{MR1022305}
for 
Sobolev-regular vector fields.

\section{Preliminaries}\label{sec:P}

\subsection{Generalities}\label{sec:gener}
We start with a brief description of the setting and  recapitulation of some of the results of \cite{MR2892584,2010arXiv1011.3273C}. Further details and references may be found in those papers (see also \cite{MR2569321,MR2283957,MR2643799} and the references therein).
In what follows, $\Rd$ is the Euclidean space of  dimension $d\ge
2$, scalar product $x\cdot y$, norm $|x|$ and Lebesgue measure $dx$.
All sets, measures and functions in $\Rd$ considered  throughout this paper will be Borel.
We denote by
$$
B(x,r)=\{y\in \Rd: |x-y|<r\},
$$ the ball of center $x\in \Rd$ and radius $r>0$.

We will consider  nonempty, bounded open set $D\subset \Rd$, whose boundary is of class $C^{1,1}$. 
The latter means that 
$r>0$ exists such that for every $Q\in \partial D$ there are balls 
$B(x',r)\subset D$ and $B(x'',r)\subset \Rd \setminus D$,
which are tangent at $Q$
(the {\it inner}\/ and {\it outer}\/  tangent ball, respectively). 
We will refer to such sets $D$ as to $C^{1,1}$ domains, without requiring connectivity. For an alternative analytic description and 
localization of $C^{1,1}$ domains we refer to \cite[Lemma~1]{MR2283957}.
We note that each connected component of $D$ contains a ball of radius $r$, and the same is true for $D^c$. Therefore $D$ and $D^c$ have a finite number of components,
which will play a role in a later discussion of extensions of the vector field to a neighborhood of $\overline{D}$. 
The distance of a given $x\in \Rd$ to $D^c$ will be denoted by $$\delta_D(x) =\inf\{|y-x|:\, y\in \Rd\setminus D\}\,.$$

{\it Constants} mean positive numbers, that do not depend on the considered arguments of the functions being compared.
Accordingly, 
notation $f(x)
\approx
g(x)$ 
means that 
there is a constant $C$ such that $C^{-1}f(x)\le g(x)\le Cf(x)$ for all $x$.
As usual, $a \land b =\min(a,b)$ and $a\vee b = \max(a,b)$. 

We will 
employ the  function space $L^2(D)$, consisting of all square integrable real valued functions, with the usual scalar product 
$$
(f,g)=\int_D f(x)g(x)dx.
$$
Generally, given $1\le p\le \infty$, the
norms in  $L^p(D)$ shall be denoted by $\|\cdot\|_{p}$. 
Customarily, $C_c^\infty(D)$ denotes the space of smooth functions on $\Rd$ with compact support in $D$.
Also, $H^1_0(D)$ 
denotes the closure of  $C_c^\infty(D)$
in the norm $|\cdot|_1$, 
where
\begin{equation}\label{dn1}
|f|_{1}^2:=\int _\Rd( |f(x)|^2+|\nabla f(x)|^2) dx
=(2\pi)^{-d}\int_\Rd (1+|\xi|^2)|\hat{f}(\xi)|^2d\xi.
\end{equation}
In the last equality we have used Plancherel theorem.  The Fourier transform of $f$ is given by
\begin{equation}
\label{fourier}
\hat{f}(\xi)=\int_\Rd e^{-i\xi\cdot x} f(x)\,dx. 
\end{equation}

\subsection{Isotropic $\alpha$-stable L\'evy process on $\Rd$}\label{sec:iso}
The discussion in Section~\ref{sec:iso} and Section~\ref{sec:kill} is valid for any $\alpha\in (0,2)$.
Let 
$$
\nu(y)
=\frac{ \alpha 2^{\alpha-1}\Gamma\big((d+\alpha)/2\big)}{\pi^{d/2}\Gamma(1-\alpha/2)}|y|^{-d-\alpha}\,,\quad y\in \Rd\,.
$$
The 
coefficient is chosen in such a way that
\begin{equation}
  \label{eq:trf}
  \int_{\Rd} \left[1-\cos(\xi\cdot y)\right]\nu(y)dy=|\xi|^\alpha\,,\quad
  \xi\in \Rd\,.
\end{equation}
We define the {\em fractional Laplacian} 
as  the $L^2(\bbR^d)$-closure of the operator 
\begin{equation*}
  \Delta^{\alpha/2}\phi(x) = 
  \lim_{\varepsilon \downarrow 0}\int_{|y|>\varepsilon}
  \left[\phi(x+y)-\phi(x)\right]\nu(y)dy\,,
  \quad
  x\in \Rd,\,\phi\in C^\infty_c(\Rd).
\end{equation*}
Its Fourier 
symbol
is given by 
$\widehat{\Delta^{\alpha/2}\phi}(\xi)=-|\xi|^{\alpha}\hat{\phi}(\xi)$, cf \eqref{fourier}.
The fractional Laplacian is 
the generator of the semigroup of the isotropic 
 $\alpha$-stable L\'evy process 
$(Y_t, \PP^x)$ on $\Rd$. Here $\PP^x$ and $\EE^x$ are the law and expectation for the process
starting at $x\in\bbR^d$. These are defined on the Borel $\sigma$-algebra of the canonical c\'adl\'ag 
{path}
space $\mathbf D([0,+\infty);\Rd)$ via transition probability densities as follows. 
We let 
$(Y_t)$ 
be
the canonical  process, i.e.
$$
Y_t(\omega):=\omega(t)\quad{\rm  for }\quad\omega \in {\mathbf D}([0,+\infty);\Rd),
$$   
and define
time-homogeneous transition 
density 
$
p(t,x,y):=p_t(y-x),
$
where 
\begin{equation}
  \label{eq:dpt}
  p_t(x)=\frac{1}{(2\pi)^d}\int_ \Rd e^{-t|\xi|^\alpha}e^{-ix\cdot\xi}\,dx\,,\quad x\in
  \Rd\,.
\end{equation}
According to 
(\ref{eq:trf}) 
and the L\'evy-Khinchine formula,
$\{p_t\}$ is a probabilistic convolution semigroup of functions with the L\'evy 
measure
$\nu(y)dy$, see e.g. 
\cite{MR2569321}.
 From (\ref{eq:dpt}) we have 
\begin{equation}
  \label{eq:sca}
  p_t(x)=t^{-d/\alpha}p_1(t^{-1/ \alpha}x)\,,\quad t>0\,,\;x\in \Rd\,.
\end{equation}
It is well-known (\cite{MR2569321}) that $p_1(x)
\approx1\land |x|^{-d-\alpha}$, hence
\begin{equation}\label{eq:oppt}
p_t(x)
\approx t^{-d/\alpha}\land \frac{t}{|x|^{d+\alpha}}
\,,\quad t>0,\,x\in \Rd\,.
\end{equation}

\subsection{
Stable process killed off $D$}\label{sec:kill}

Let $\tau_D=\inf\{t>0: \, Y_t\notin D\}\,$ 
be the {\it time of the first exit}\/ of the (canonical) process from $D$.
For each $(t,x)\in(0,+\infty)\times D$, the measure 
$$
A\mapsto P_D(t,x,A):=\PP^x[Y_t\in A,\tau_D>t]
$$ is absolutely continuous 
with respect to the Lebesgue measure
on $D$. Its density $p_D(t,x,y)$ is continuous in $(t,x,y)\in(0,+\infty)\times D^2$ and 
satisfies  G.~Hunt's formula (see \cite{MR1490808}, \cite{MR2283957}),
\begin{equation}\label{eq:Hunt:1}
p_D(t,x,y) = p(t,x,y) - \bbE^x\left[\tau_D < t;\; p(t-\tau_D,
  Y_{\tau_D},y)\right].
\end{equation}
In addition, the kernel is symmetric (see \cite{MR1671973, MR2283957} for discussion and references):
\begin{equation}
\label{symmetric}
p_D(t,x,y) =p_D(t,y,x) ,\quad
t>0,\,x,y\in D.
\end{equation}
and  defines a strongly continuous semigroup on $L^2(D)$,
$$
P_t^Df(x):=\int_{D}p_D(t,x,y)f(y)dy.
$$
We 
shall denote $P_t=P_t^\Rd$.
The {\em Green function} of $\Delta^{\alpha/2}$ for $D$ is defined as
$$
G_D(x,y)=\int_0^\infty 
p_D(t,x,y)dt, 
$$
and the respective {\em Green operator} is
\begin{eqnarray*}
G_Df(x) &:=& \EE^x\left[ \int_0^{\tau_D}f(Y_t)dt\right]
=\int_{\Rd}G_D(x,y)f(y)dy.
\end{eqnarray*}
We have (\cite{MR2892584}),
\begin{equation}\label{eq:fg}
G_D (\Delta^{\alpha/2}f)=-f\,,\quad  f \in C^\infty_c(D),
\end{equation}
and 
\begin{equation}\label{eq:fg1}
\Delta^{\alpha/2}(G_D f)=-f,\quad f \in L^2(D),
\end{equation}
see \cite[Lemma~5.3]{MR1825645}, where the fractional laplacian
operator  appearing above is defined in the sense of
distributions theory. See formula \eqref{left-inv} below for another
statement.

Thus the generator of $(P^D_t)$ is 
$\Delta^{\al/2}$ with zero (Dirichlet) exterior conditions, i.e. with
the domain equal to the range of $G_D$,
\begin{equation}
\label{domain}
{\mathcal D}(\Delta^{\alpha/2})=G_D(L^2(D))=:R(G_D). 
\end{equation}
The following estimate 
{has been} proved by Kulczycki \cite{MR1490808} and Chen and
Song \cite{MR1654824} (see also \cite[Theorem 21]{MR1991120}),

\begin{eqnarray}
 \label{GreenEstimates}
  &&G_D(x,y)\, 
\approx\;  |x-y|^{\alpha-d}
\frac{\delta_D(x)^{\alpha/2}\delta_D(y)^{\alpha/2}}{[\delta_D(x)\vee |x-y| \vee \delta_D(y)]^\alpha}\label{GreenEstimatesTJ}\\
&&  \approx\;  |x-y|^{\alpha-d}\left(\frac{\delta_D(x)^{\alpha/2}\delta_D(y)^{\alpha/2}}{|x-y|^\alpha} \land 1\right)\nonumber
\,, \qquad x,y \in D\,.\nonumber
\end{eqnarray}
In particular,
$$
G_D(x,y)\leq C |x-y|^{\alpha-d},\, \qquad x,y \in D.
$$
From \eqref{GreenEstimates}
it also follows that 
\begin{align}\label{oco}
\EE_x \tau_D =G_D 1(x)\approx \delta_D^{\alpha/2}(x), \quad x\in D.
\end{align}
{From \cite[Corollary 3.3]{MR1936936} we have the following gradient estimate,
\begin{equation}\label{eq:GradEstimGreen}
|\nabla_x G_D(x,y)| \le d \frac{G_D(x,y)}{\delta_D(x) \land |x-y|}\,,
\quad x,y \in D,\,\; x\neq y\,.
\end{equation}}
We define $H^{\alpha/2}=H^{\alpha/2}(\Rd)$ as the subspace of $ L^2(\Rd)$ made of
those elements for which
\begin{equation}
\label{012210}
|f|_{\alpha/2}:=\{\|f\|_2^2+{\cal E}^{\alpha}(f,f)\}^{1/2}<+\infty.
\end{equation}
Here (the Dirichlet form) ${\cal E}^{\alpha}$ is given as follows (cf \cite{MR2778606}):
\begin{eqnarray}
\nonumber
&&
{\cal E}^{\alpha}(f,g):=\lim_{t\to 0+}\frac{1}{t}((I-P_t)f,f)\\
&&
=\frac12 \int_\Rd \int_\Rd [f(x)-f(y)]^2\nu(y-x)dxdy,\label{011109b}
\end{eqnarray}
The above formula can be also used to define a bilinear form ${\cal E}_{\al}(\cdot,\cdot)$ 
on $H^{\al/2}\times H^{\al/2}$ by polarization. We also have
\begin{equation}\label{eq:GfDa}
|f|_{\alpha/2}^2=(2\pi)^{-d}\int_\Rd
(1+|\xi|^\alpha)|\hat{f}(\xi)|^2d\xi,
\end{equation}
therefore
\begin{equation}\label{eq:ca1}
|f|_{\alpha/2}\le \sqrt{2}|f|_1,
\quad f\in L^2(\Rd).
\end{equation}
We define $H^{\alpha/2}_0(D)$ as the closure of $C_c^\infty(D)$ in the norm 
$|\cdot|_{\al/2}$.
By Theorems 4.4.2,  A.2.10 and formula (4.3.1) of 
\cite{MR2778606},
we have (cf. \eqref{011109b})),
$$
{\cal E}^{\alpha}(f,
g)=\lim_{t\to0+}\frac{1}{t}((I-P_t^D)f,
g),\quad 
 f,g\in H^{\alpha/2}_0(D).
$$
If $f\in H_0^{\alpha/2}(D)$ and $g$ belongs to the domain of the
fractional Laplacian,
then
\begin{equation}\label{eq:GfD}
{\cal E}^\alpha(f,g)=-\int_\Rd f(x)\Delta^{\alpha/2}g(x)dx .
\end{equation}

\subsection{Gradient perturbations of $\Delta^{\alpha/2}$ on $\Rd$}
Throughout the remainder of the paper we always assume that 
$1<\alpha<2$, and $b:\Rd\to \Rd$ is a bounded vector field.
For $t>0$ and $x,y \in \RR^d$ we let
$$
 p_0(t,x,y)  :=  p(t,x,y)
 $$
 and for each $n\ge1$
\begin{equation}\label{eq:dipn}
 p_n(t,x,y)  :=  \int_0^t \int_{\RR^d} p_{n-1}(t-s,x,z) b(z) \cdot \nabla_z p(s,z,y)\,dz\,ds.
\end{equation}

Let
\begin{equation}
\label{010806}
\tp(t,x,y)=\sum_{n=0}^\infty p_n(t,x,y).
\end{equation}
It follows from  \cite[Theorem~2 and Example~2]{MR2643799} that 
series \eqref{010806}
converges   uniformly on compact subsets of $(0,+\infty)\times (\bbR^d)^2$.
From the results of  
\cite{MR2283957}, we know that $\tp(t,x,y)$ is a transition probability density function, i.e. it is non-negative and
\begin{equation}
\label{030806}
\int_{\bbR^d}\tilde p(t,x,y)dy=1,\quad
t>0,\,x\in \RR^d,
\end{equation}
In addition, $\tilde p$ is continuous  on $(0,+\infty)\times (\bbR^d)^2$, and
\begin{equation}\label{ptxy_comp}
  c_T^{-1} p(t,x,y) \le \tp(t,x,y) \le c_Tp(t,x,y)\,,\qquad x,y \in \RR^d\,,\; 0<t\le T\,,
\end{equation}
where $c_T \to 1$ if $T \to 0$.
In fact, this holds under much weaker, Kato-type condition on $b$, see \cite[Theorems~1 and 2]{MR2283957}. 

We denote by $\tPP^x$ and $\tEE^x$ the law  and expectation on  ${\mathbf D}([0,+\infty);\Rd)$ for the (canonical) Markov process starting at $x$
and defined by the transition probability density $\tp$,

\begin{remark}
\rm
$\tPP^x$ may also be defined by solving
stochastic differential equation $dX_t=dY_t+b(X_t)dt$. Such equations
have been studied in dimension $1$ in \cite{MR0368146} under the assumptions of
boundedness and continuity of the vector field;
also for
$\alpha=1$. We  refer the reader to \cite[formula
(13)]{MR1310558}, for a closer description of a  connection to \eqref{010806} and \eqref{eq:dipn}.
\end{remark}
\begin{remark}
\rm
We may also define the perturbation series for $-b(x)$.
In what follows, objects pertaining to $-b$ will be marked with the superscript hash ($^\#$),
e.g. 
$
\tp^\#(t,x,y)=\sum_{n=0}^\infty (-1)^n p_n(t,x,y).
$
\end{remark}

\subsection{Antisymmetry of the perturbation}
If $\diver b=0$ on $D$ in the sense of distributions theory, and $f\in C_c^\infty(D)$, then
\begin{equation}\label{egas}
\int_D f(x)b(x)\cdot \nabla f(x)dx=\frac12\int_D b(x)\cdot \nabla f^2(x)dx=0.
\end{equation}
If also $g\in C^\infty_c(D)$ and {we substitute  $f+g$ for $f$ in \eqref{egas},  then we obtain} 
\begin{equation}\label{antis}
\int_D f(x)b(x)\cdot \nabla g(x)dx=-\int_D g(x)b(x)\cdot \nabla f(x)dx.
\end{equation}
The last equality extends to 
arbitrary $f,g\in H_0^1(D)$. 

\begin{prop}\label{lem:duality}
If $\diver b=0$ on $\Rd$, then we have
$p^{\#}_n(t,x,y)=p_n(t,y,x)$
and $\tp^{\#}(t,x,y)=\tp(t,y,x)$
for all $t>0$, $x,y\in \Rd$, and $n\ge0$.
\end{prop}
\proof
Let $s,t>0$, $x,y\in \Rd$. By \eqref{antis} (see also \cite[Lemma~4]{MR2875353}), 
\begin{equation}\label{eq:ip}
\int_\Rd p(t,x,z)b(z)\cdot\nabla_z p(s,z,y)dz=
-\int_\Rd p(s,y,z)b(z)\cdot\nabla_z p(t,z,x).
\end{equation}
From \eqref{eq:ip} we conclude that $p_1^{\#}(t,x,y)=p_1(t,y,x)$.
For a general $n \ge 1$, by 
\eqref{eq:dipn} we have
\begin{equation}\label{eq:p_n}
p_n(t,x,y) = \int_{D_n(t) }\ d{\bf s}\int_{(\RR^d)^n} P_n({\bf s},x,{\bf z},y)\;
d{\bf z}.
\end{equation}
where ${\bf z}:=(z_1,\ldots,z_n)$, ${\bf s}:=(s_1,\ldots,s_n)$, 
$$
D_n(t):=[{\bf s}:0\le s_1\le \ldots\le s_n\le t],
$$
and
\begin{equation}\label{def:p_n_ext}
P_n({\bf s},x,{\bf z},y)
= p(s_1,x,z_1)\prod_{i=1}^{n}\left[b(z_i) \cdot \nabla_{z_i} p(s_{i+1}-s_i,z_i,z_{i+1})\right],
\end{equation}
with the convention that $(s_{n+1},z_{n+1})=(t,y)$, cf. \cite{(4.6)]{MR2875353}}.
Using formula \eqref{eq:ip} $n$ times in space, and then integrating in time we see that 
\begin{equation}\label{eq:altpn}
p_n(t,x,y)=(-1)^n p_n(t,y,x),
\end{equation}
which yields the 
identities stated  in the lemma.
\qed

\begin{remark}
\rm
A strengthening of Proposition~\ref{lem:duality} will be given in Corollary~\ref{cor010908} below.
\end{remark}
\begin{remark}
\rm
If $\diver b=0$ on $\Rd$, then
\eqref{030806} 
and 
Proposition  \ref{lem:duality} 
yield
\begin{equation}
\label{020806}
\int_{\bbR^d}\tilde p(t,x,y)dx=1,\quad
t>0,\,y\in \RR^d.
\end{equation}
\end{remark}

\subsection{Gradient perturbations with Dirichlet conditions}\label{sec:deftG}
We recall that $D$ is a bounded $C^{1,1}$ domain in $\Rd$. Hunt's
  formula 
may 
be used to define the transition
probability density of the (first) perturbed and (then) killed process
\cite{MR2892584}. 
Thus, for $t>0$, $x,y\in D^2$ we let
\begin{equation}\label{eq:Hunt}
\tp_D(t,x,y) = \tp(t,x,y) - \tEE^x\left[\tau_D < t;\; \tp(t-\tau_D,
  Y_{\tau_D},y)\right]\,.
\end{equation}
We have 
$$
\tEE^x[\tau_D < t; f(Y_t)]=\int_D \tp_D(t,x,y)f(y)dy,\quad
t>0,\,x\in D,\,f\in L^\infty(D),
$$
and
\begin{equation}\label{eq:eotp1}
0\le \tp_D(t,x,y)\le \tp(t,x,y)\le c_t p(t,x,y),\quad
t>0,\,x,y\in D,
\end{equation}
where 
{$\sup_{t\in[0,T]}
c_t<+\infty$  for any $T>0$. }

By \cite[formula (40)]{MR2892584}, there exist constants $c,C>0$ such that
\begin{equation}\label{eq:eotp}
\tp_D(t,x,y) \le Ce^{-ct}\,,\quad t\ge 1, \quad x,y \in
  \RR^d.
\end{equation}

We define the Green function 
of 
$D$ 
for $L$:
\begin{equation}\label{eq:deftG}
\tG_D(x,y)=\int_0^\infty \tp_D(t,x,y)dt\,.
\end{equation}
Clearly $\tG_D$
is nonnegative.
By Blumenthal's 0-1 law, $\tp_D(t,x,y)=0$ for all $t>0$ (and thus $\tG_D(x,y)=0$) if $x\in D^c$ or $y\in D^c$, see (\ref{ptxy_comp}).
It follows from \eqref{eq:egf} and \eqref{GreenEstimatesTJ} (see  \cite[Lemma 7]{MR2892584}), that
\begin{equation}\label{eq:ubRk}
\tG_D(x,y) \le C_0|x-y|^{\alpha-d}\,,\quad x,y\in D.
\end{equation}
The main result of \cite{MR2892584}  asserts that 
  \begin{equation}
    \label{eq:egf}
\CIII^{-1}G_D(x,y) \le \tilde G_D(x,y) \le \CIII G_D(x,y),\quad \mbox{ if } x,y \in D
\,,
  \end{equation}
and $\tG_D(x,y)$ is continuous for $x\neq y$ {(for estimates of $\tp_D$ see \cite{2010arXiv1011.3273C})}.
We consider the 
integral operators
$$
\tP^{D}_t f(x)=\int_\Rd \tp_D(t,x,y)f(y)dy\,,\quad t>0,
$$
and 
\begin{equation}\label{eq:doG}
\tG_D f(x)=\int_\Rd \tG_D(x,y)f(y)dy. 
\end{equation}
In light of \eqref{eq:eotp1} and \eqref{eq:egf}, the above 
operators are bounded on every $L^p(D)$, $p\in[1,+\infty)$. The analogous
operators $\tilde P^{D{\#}}_t$ and $\tG^{\#}_D$, defined for $-b(x)$, turn out to be
mutually adjoint on $L^2(D)$,
as {follows from}
Corollary \ref{cor010908} below.
From \eqref{030806}, \eqref{020806} and \eqref{eq:Hunt}
we obtain that
$\{\tP^{D}_t \}$ is a semigroup of contractions in 
every $L^p(D)$, $p\in[1,+\infty]$. It is strongly continuous for $p\in[1,+\infty)$.
Let $L$ be the $L^2(D)$ generator of the semigroup, with the domain ${\mathcal D}(L)=R(\tG_D)=\tG_D(L^2(D))$. We have 
\begin{equation}
\label{left-inv}
L\tilde G_Df=-f,\quad
f\in L^2(D).
\end{equation}
It has been shown in \cite{MR2892584}  that the following crucial  recursive formula
holds
\begin{equation}\label{eq:pf}
  \tG_D = G_D + \tG_D (b\cdot\nabla) G_D,
\end{equation}
and for all $\varphi \in C_c^\infty(D)$ and $x\in D$ we have 
\begin{equation}\label{tG_D_fs}
\int_{D} \tG_D(x,z) \left(\Delta^{\alpha/2}\varphi(z) + b(z) \cdot \nabla \varphi(z) \right)\,dz = - \varphi(x)\,.
\end{equation}
Later on we shall also consider 
the operator
\begin{equation}\label{eq:defL:a}
L_A = \Delta^{\alpha/2} + A b\cdot\nabla\,,
\end{equation}
corresponding to the vector field $Ab$, where $A\in\bbR$ and we
  let $A\to \infty$. 
Clearly, if $A$ is fixed, then there is no loss of generality to focus on $L=L_1$.

\section{Comparison of the domains of generators}\label{sec:cod}

\label{sec4}

The following pointwise version of \eqref{eq:pf} is proved in \cite[Lemma 12]{MR2892584},
\begin{equation}
\label{011307}
\tG_D(x,y)=G_D(x,y)+\int_D\tG_D(x,z)H(z,y)dz,\quad
x,y\in D,
\end{equation} 
with
\begin{equation}\label{defHxy}
H(x,y):= b(x)\cdot \nabla_x G_D(x,y), \quad x,y\in D.
\end{equation}
We define
\begin{equation}\label{eq:defH}
Hf(x)=\int_DH(x,y)f(y)dy.
\end{equation}
After a series of auxiliary estimates, we will prove $H$ to be compact on $L^2(D)$.
\begin{lem}
\label{lem011107}
There exist $C>0$ such that
 \begin{equation}
\label{041107}
|H(x,y)|\le C\delta_D^{\al/2-1}(x) |y-x|^{\al-1-d}\delta_D^{1-\al/2}(y),\quad 
x,y\in D.
\end{equation}
\end{lem}
\proof
Using (\ref{eq:GradEstimGreen}) we obtain that
$$
|H(x,y)|\le c \dfrac{G_D(x,y)}{\delta_D(x) \land |x-y|}=c\dfrac{G_D(x,y)[\delta_D(x)\vee |x-y|]}{\delta_D(x)|x-y|},
$$
which we bound from above, 
thanks to (\ref{GreenEstimatesTJ}), by
\begin{align}
\label{051107}
&c_1
\delta_D(x)^{\alpha/2-1}|x-y|^{\alpha-1-d}
\dfrac{[\delta_D(x)\vee |x-y|]\delta_D(y)^{\alpha/2}}{[\delta_D(x)\vee
  |x-y| \vee \delta_D(y)]^{\alpha}}\,   \nonumber\\
&\le c_1\,\delta_D(x)^{\alpha/2-1}|x-y|^{\alpha-1-d}
\dfrac{\delta_D(y)^{\alpha/2}}{[\delta_D(x)\vee
  |x-y| \vee \delta_D(y)]^{\alpha-1}}\,,  
\end{align}
and this yields \eqref{041107}.
\qed

\begin{lem}\label{lem051307}
If $0<1-\beta<\gamma<1$,
then there is $c=c(d,\beta,\gamma,D)$ such that
\begin{equation}
\label{051307}
\int_D|y-x|^{-(d-1)-\beta}\delta_D(y)^{-\gamma}dy\le c \delta_D(x)^{1-\beta-\gamma}, \quad
x\in D.
\end{equation}
\end{lem}
\proof
We only need to examine points  close to $\partial D$. { Given such a
point   we consider the
  integral over its neighborhood ${\cal O}$.
The neighborhood can be chosen in such a way that, after a
bi-Lipschitz change of variables (see \cite[formula
(75)]{MR2892584}, \cite[formula (11)]{MR1936936}), we can reduce our
consideration to the case when $\partial D \cap \overline{{\cal O}}
\subset [x_1=0]$, ${\cal O}\subset [x_1>0]$ and $\delta_D(x)=x_1$
for $x\in {\cal O}$. }
The respective
integral over {${\cal O}$ }
is then estimated by
\begin{equation*}
\label{061307}
c\int_0^\infty 
\frac{dy_1}{y_1^{\gamma}}\int_{\RR^{d-1}}
\frac{dx'}
{\left(|x_1-y_1|+|x'|\right)^{d-1+\beta}}
=
c_1\int_0^\infty |y_1-x_1|^{-\beta}y_1^{-\gamma}dy_1=c x_1^{1-\beta-\gamma}.\qquad\qed
\end{equation*}
\begin{remark}
{\rm
The result is valid for all bounded open {\it Lipschitz} sets (\cite{MR1936936}) in all dimensions $d\in \NN$.
}
\end{remark}
Lemma~\ref{lem051307} yields the following.
\begin{cor}\label{cordk}
Let $1<\alpha'\le \alpha<2$ and $\al'-1<\jo,\ka<1$. Suppose also that
\begin{align}
\label{010312}
K(x,y)&:=\delta_D(x)^{\al/2-1}|y-x|^{-(d-1)-(2-\al')}\delta_D(y)^{1-\al/2},\\
p(x)&:=\de(x)^{1-\al/2-\jo}\quad \mbox{ and } \quad  q(y):=\de(y)^{\al/2-1-\ka},\qquad x,y\in D.\nonumber
\end{align}
Denote 
$$
pK(y)=\int_D p(x)K(x,y)dx, \quad Kq(x)=\int_D K(x,y)q(y)dy.
$$
Then,
\begin{equation*}
pK(y)\le c_1 \de_D(y)^{\al'-\al-\jo},\quad 
Kq(x)\le c_2 \de_D(x)^{-2+\al'+\al/2-\ka},\quad x,y\in D.
\end{equation*} 
{Constants  $c_1, c_2$ depend only on $\al,\al',\jo,d,D.$}
\end{cor}
{
The above result shall be used to establish that operator $K$ (thus
also $H$) is $L^2$ bounded via Schur's test, see \cite[Theorem 5.2]{MR517709}.
Note that we have
 $pK\le cq$, for some constant $c>0$, provided that
$\jo-\ka\le \al'-\al+1$. Likewise,  $Kq\le c p$, provided that $\jo-\ka\ge 3-\al'-\al$.
Summarizing, we will require the following conditions:
\begin{align}
&1<\al'\le \al<2,\label{C1}\\
&\al'-1<\ka,\jo<1,\quad \mbox{and}\label{C2}\\
&3-\al'-\al \le \jo-\ka \le \al'-\al +1.\label{C3}
\end{align}
\begin{lem}\label{lemary}
Conditions \rm (\ref{C1})-(\ref{C3})
hold, 
if $1<\al<2$ and
$$
\al' :=\al-\ga(\al),\quad
\ka :=\al-1+\ga(\al),\quad
\jo :=1-\ga(\al),
$$
where $\ga(\al):=(\al-1)(2-\al)/3$. 
\end{lem}
\proof
We have
$
\jo-\ka-(3-\al'-\al)=
1-\al(2-\al)\ge 0
$, etc.
\qed
\begin{prop}
\label{prop021307}
Operator $H$ is compact on $L^2(D)$.
\end{prop}
\proof
We note that $|H(x,y)|\le c K(x,y)$, where $K(x,y)$ is given by
\eqref{010312} and $c>0$ is some constant.
Let $\al',\ka,\jo$ be  as in Lemma~\ref{lemary}.
By  Corollary~\ref{cordk} and the discussion preceding Lemma~\ref{lemary},
there are constants $c_1(\al'), c_2(\al')$ such that
\begin{equation}
\label{050312}
pK\le c_1(\al')q\quad\mbox{ and }\quad Kq\le c_2(\al')p
\end{equation}
Using Schur's test for boundedness of integral operators on $L^2$
  spaces (see \cite[Theorem 5.2]{MR517709}) we obtain 
$$\|H\|\le c\|K\|\le c \sqrt{c_1(\al')c_2(\al')}<\infty.$$
For $N>0$, $r\in \RR$ we let $\Phi_N(r)=(r\wedge N)\vee (-N)$, $r\in\bbR$, and
$$
H_N(x,y):=\de_D^{\al/2-1}(x)\Phi_N\left(\de_D(x)^{1-\al/2}H(x,y)\de_D(y)^{\al/2-1}\right) \de_D^{1-\al/2}(y).
$$
Each operator  $H_Nf(x):=\int_D H_N(x,y)f(y)dy$ is Hilbert-Schmidt,
hence compact (\cite[Theorem 4, p. 247]{MR1892228}), due to the fact that
$$
\int_{D^2}H_N^2(x,y)dxdy\le
N^2\int_D\de_D(x)^{\al-2}dx\int_D\de_D(y)^{2-\al}dy
<\infty.
$$
We  have $H(x,y)=H_N(x,y)$, unless
$\de_D(x)^{1-\al/2}|H(x,y)|\de_D(y)^{\al/2-1}>N$.
The latter may only happen if $x,y$ satisfy
$
C|y-x|^{-(d+1-\al)}>N$, with $C>0$ 
(cf \eqref{041107}), or
equivalently when
$
|y-x|<(C/N)^{1/(d+1-\al)}
$,
and
then there are constants $C,C_1>0$ such that 
\begin{align*}
&|H(x,y)-H_N(x,y)|\le |H(x,y)|\le C |x-y|^{\al'-\al}K(x,y)\\
&\le
C_1N^{(\al'-\al)/(d+1-\al)} K(x,y).
\end{align*}
}
This estimate for $H(x,y)-H_N(x,y)$
actually holds for all $x,y\in D$, hence
\eqref{050312} yields
\begin{align*}
p|H-H_N|&\le C_1N^{(\al'-\al)/(d+1-\al)}c_1(\al')q,\\
|H-H_N|q&\le C_1N^{(\al'-\al)/(d+1-\al)}c_2(\al')p.
\end{align*}
Applying Schur's test we get $\|H-H_N\|\le C_1N^{(\al'-\al)/(d+1-\al)}\sqrt{c_1(\alpha')c_2(\alpha')}$.
Since $H$ is a norm limit of compact operators, it is compact, 
too.
\qed

In view of Proposition\ref{prop021307}
we may regard   the gradient 
operator $b\cdot \nabla$ 
as a small perturbation of $\Delta^{\alpha/2}$ when $\al\in(1,2)$. In fact $b\cdot \nabla$ is relatively compact with respect to  $\Delta^{\alpha/2}$ with Dirichlet conditions in the sense of \cite[IV.1.3]{MR1335452}.

\begin{thm}
\label{thm011307}
${\mathcal D}(L)={\mathcal D}(\Delta^{\al/2})$.
\end{thm}
\proof
By virtue of \eqref{011307}, we can write
\begin{equation}
G_D=\tG_D(I-H)\quad \mbox{on }L^2(D).\label{eq:30212}
\end{equation}
By \eqref{eq:fg1}, $G_D$ is injective, therefore so is $I-H$.
In particular, $1$ is not an eigenvalue of $H$.
Since $H$ is compact, by 
the Riesz-Schauder theory \cite[Theorem X.5.1, p. 283]{MR617913},
$I-H$ is invertible.
Thus,
$$
{\mathcal D}(\Delta^{\alpha/2})=R(G_D)=R(\tG_D(I-H))=R(\tG_D)={\mathcal D}(L). 
$$
For future reference we remark that 
$I-H^*$ 
is invertible, too.
\qed

\begin{cor}
\label{cor011307}
$ {\mathcal D}(L)\subset H^1_0(D)$.
\end{cor}
\proof 
By \cite[Lemma 10]{MR2892584},
\begin{equation}\label{eq:cdiG}
\nabla G_Dg(x)=\int_D\nabla_x G_D(x,y)g(y)dy,
\end{equation}
for any bounded function $g$.
Invoking the argument used in Proposition \ref{prop021307}, we conclude that there is a number $c$ independent of $g$ for which
\begin{equation}\label{eq:odnp}
\|\nabla G_Dg\|_2\le c\|g\|_2.
\end{equation}
By approximation, 
 \eqref{eq:cdiG} and \eqref{eq:odnp} extend to all $g\in L^2(D)$. Furthermore,
if $\{g_j\}\subset C_c^\infty(D)$ and $\lim_{j\to+\infty}g_j= g$ in $L^2(D)$, then $\lim_{j\to+\infty}
|G_Dg-G_Dg_j|_1= 0$, cf. \eqref{dn1}.
Each $G_Dg_j$ is continuous on $\overline D$ and $G_Dg_j=0$ on $\partial D$, hence $G_Dg_j\in H_0^1(D)$ (see \cite[Theorem 5.37]{MR2424078}).
Therefore, 
$
G_Dg\in H_0^1(D)$.
\qed

The following result justifies our notation $L=\Delta^{\alpha/2}+b\cdot \nabla$.
\begin{prop}
\label{prop011009}
If $f \in {\mathcal D}(L)$, then
$Lf=\Delta^{\al/2} f+b\cdot\nabla f$.
\end{prop}
\proof
Let $f\in {\mathcal D}(L)$ and  $h=b\cdot\nabla f$.
Applying $\tilde G_D$  to $\Delta^{\al/2} f+h$ 
and 
using \eqref{eq:pf}, we obtain
\begin{eqnarray}
\label{021009}
&&
\tilde G_D(\Delta^{\al/2} f+h)=( G_D+\tilde G_D(b\cdot \nabla)G_D)(\Delta^{\al/2} f+h)\\
&&
=-f+G_Dh-\tilde G_Dh+\tilde G_D(b\cdot \nabla)G_D h=-f,\nonumber
\end{eqnarray}
which, thanks to \eqref{left-inv}, concludes the proof.
\qed

We shall also observe
the following localization principle 
{for}
our perturbation problem.
\begin{cor}
\label{cor010212}
Operators $\tilde G_D$ and $(\tP_t^D, t>0)$ do not depend 
on the values of $b$ on $D^c$.
\end{cor}
\proof
If $b_1$ and $b_2$ are two (bounded) vector fields on $\Rd$ equal to
$b$  
on $D$,
and $\tilde G_D^{(1)}$ and $\tilde G_D^{(2)}$ are the  corresponding Green's operators,
then, {by \eqref{eq:30212}},
$$
\tilde G_D^{(1)}(I-H)=\tilde G_D^{(2)}(I-H),
$$
which identifies $\tG_D:=\tG_D^{(1)}=
\tG_D^{(2)}$ as operators on $L^2(D)$, hence also as 
functions defined in Section~\ref{sec:deftG}.
A similar 
conclusion 
for 
$\tP_t^D$ follows 
from the fact that $\tilde G_D^{-1}$ is the generator of the 
semigroup.
\qed

We will identify the adjoint operator of $H$ on $L^2(D)$.
\begin{lem}\label{lemHs}
$H^*=-G_Db\cdot \nabla$.
\end{lem}
\proof
If $f,g\in L^2(D)$, then by Corollary~\ref{cor011307} and \eqref{antis},
\begin{align*}
&{(Hf,G_D g )}=\left( b\cdot \nabla G_Df, G_D g \right)=-\left(b\cdot \nabla G_D g, G_D f\right)
=-\left(f,G_D b\cdot \nabla G_D g \right).
\end{align*}
Functions $\{G_D g,\,g\in  L^2(D)\}$ form a dense set in $L^2(D)$, which ends the proof.
\qed

When defining $\tP^D_t$, via \eqref{eq:Hunt},  
we may encounter the situation when  $b$ is given only on $D$. 
We may extend  the
field outside $D$ by letting e.g.  $b=0$ on $D^c$. Of course such an
extension  needs not  satisfy ${\rm div b}=0$, even though the condition may
hold on $D$. However, we still 
have a
local analogue of Proposition \ref{lem:duality}. Recall that   the transition probability densities and Green function corresponding to the vector field $-b$ are marked with a hash ($^\#$).
\begin{cor}
\label{cor010908}
If  ${\rm div}\,b
=0$ on $D$, then
for all $t>0$ and $x,y\in D$,
\begin{equation}
\label{020908}
 \tp^{\#}_D(t,x,y)=\tp_D(t,y,x),
\end{equation}
 and 
 \begin{equation}
 \label{010707a}
 \tG^{\#}_D(x,y)=\tG_D(y,x).
 \end{equation}
\end{cor}
\proof 
We shall first prove \eqref{010707a}, or, equivalently, that
 \begin{equation}
 \label{010707b}
 \tG^{\#}_D=\tG^*_D,
 \end{equation}
where 
$\tG^*_D$ is 
adjoint to  $ \tG_D$ 
on $L^2(D)$.
 From\eqref{left-inv} applied to $ \tG^{\#}_D$  we 
have,
$$
 (\Delta^{\al/2}-b\cdot \nabla) \tG^{\#}_D=-I.
$$
Applying the operator $G_D$  
from the left 
to both sides of the 
equality we obtain
\begin{equation}
  \label{eq:010212}
  \tG^{\#}_D+G_D (b\cdot \nabla) \tG^{\#}_D=G_D.
\end{equation}
By Lemma~\ref{lemHs},
\begin{equation}
  \label{eq:010212a}
  G_D=(I-H^*)\tG_D^{\#}.
\end{equation}
Taking adjoints of both sides of \eqref{eq:30212} we also obtain,
\begin{equation}
  \label{eq:020212}
   G_D=(I-H^*)G^*_D.
\end{equation}
We have already noted in the proof of Theorem \ref{thm011307}
that $I-H^*$ is 
a linear automorphism of $L^2(D)$. Therefore 
\eqref{eq:010212a} and \eqref{eq:020212} 
give \eqref{010707b}. 
Furthermore, let
$L^*$ denote the adjoint of $L$ on $L^2(D)$.
We 
have
\begin{equation}
\label{010512}
L^{\#}=(\tG_D^{\#})^{-1}=(G^*_D)^{-1}=L^{*},
\end{equation}
where the last equality follows from \cite[Lemma XII.1.6]{MR0188745}.
Let $(\tP_t^{D{*}})$ be the semigroup adjoint to $(\tilde P_t^D)$.
By \cite[Corollary 4.3.7]{MR1784426}, the generator of $(\tP_t^{D{*}})$ is $L^{*}$.
Since $L^{\#}=L^{*}$, the semigroups are equal. The corresponding  kernels are defined pointwise, therefore they satisfy \eqref{020908}.
\qed

\section{
Krein-Rutman eigen-pair}
\label{sec3}

\begin{lem}\label{lem:GDc1}
Operator $\tilde P_t^D$ is  Hilbert-Schmidt  on $L^2(D)$ for each $t>0$. 
\end{lem}
\proof
Let $\{e_n,\,n\ge0\}$ 
be an orthonormal base in $L^2(D)$. 
By \eqref{eq:eotp1} and 
Plancherel's identity, we bound the Hilbert-Schmidt norm as follows.
Using \eqref{ptxy_comp} and \eqref{eq:oppt}
\begin{align*}
&\sum_{n=0}^{+\infty}\|\tilde P_t^De_n\|_2^2 \nonumber
=\int_{D} \sum_{n=0}^{+\infty}\left(\int_{D}\tilde p_D(t,x,y)e_n(y)dy\right)^2dx\\
&=\int_{D} dx\int_{D}\tilde p_D^2(t,x,y)dy
\le \int_{\Rd} dx\int_{D}
\tp^2(t,x,y)dy
\nonumber\\
&
\le c_t^2|D|\int_{\Rd}p_t^2(x)dx<+\infty.\qquad\qed
\label{eq:onL2}
\end{align*}

\begin{lem}\label{lem:GDc1KB}
$\tG_D$ is compact on $L^2(D)$. 
\end{lem}
\proof
Let $N>0$, and
$\tilde G_D^{(N)}(x,y)=\tilde G_D(x,y)\wedge N$.
The integral operator $\tG_D^{(N)}$ on $L^2(D)$ corresponding to the kernel $\tilde G_D^{(N)}(x,y)
$ is compact. Indeed,  it 
has a finite
Hilbert-Schmidt norm: 
$$
\int_D\int_D [\tilde G_D^{(N)}(x,y)]^2dxdy\le N^2|D|^2<\infty.
$$
The  norm of 
$
\tG_D-\tG_D^{(N)}
$ 
on $L^2(D)$ may be directly estimated as follows, 
\begin{eqnarray}
\label{012208}
&&
\|\tG_D-\tG_D^{(N)}\|^2
\le \left\{\sup_x \int_D\left|  \tG_D(x,y) - \tilde G_D^{(N)}(x,y)\right|dy\right\} \\
&&\times \left\{\sup_y \int_D \left| \tG_D(x,y) - \tilde G_D^{(N)}(x,y)\right|dx\right\},\nonumber
\end{eqnarray}
see, e.g., Theorem 3, p. 176 of \cite{MR1892228}.
We let $N\to \infty$.  The functions $\tG_D(x,\cdot)$ (and $\tG_D(\cdot,y)$) are uniformly integrable on $D$ by \eqref{eq:ubRk}.
Since $\tG_D$ is approximated in 
the norm topology by compact operators,  
it is compact.
\qed

In the 
special case when  $b(x)\equiv0$, 
$(\tP_t^D)$ equals $(P_t^D)$,
a symmetric contraction semigroup on $L^2(D)$, whence
the Green operator $G_D$ is symmetric, compact and positive definite.  
The spectral theorem 
yields the following.
\begin{cor}
\label{cor011206}
$G_D^{\beta}$ is symmetric and compact for every $\beta>0$.
\end{cor}

 By \eqref{eq:egf} and \eqref{GreenEstimatesTJ}, $\tG_D$ is also irreducible.
 Krein-Rutman theorem (see \cite{MR0038008}) implies that there exists a unique nonnegative $\phi\in \Lt$ and a number $\lambda>0$ such that $\|\phi\|_2=1$ and 
 \begin{equation}\label{eq:Gein}
\tG_D\phi=\frac1\lambda \phi\,,
\end{equation}
We shall call 
$(\lambda,\phi)$  the {\em principal eigenpair} corresponding to  $L$.
From \eqref{left-inv} we have
$\phi\in {\mathcal D}(L)$ and $L\phi=-\lambda\phi$. Furthermore,
$$
\tP_t^D\phi=e^{-\lambda t}\phi,\quad 
t\ge0.
$$ 
The formula \eqref{eq:doG} yields extra regularity of $\phi$, 
as follows.
\begin{lem}
\label{lem011106} 
If \eqref{eq:Gein} holds for some $\phi\in L^1(D)$ and $\lambda\not=0$,
then $\phi\in C(\bar D)$ and there is $C=C(\alpha,b,D,\lambda)$ such that
\begin{equation}
\label{011106}
|\phi(x)|\le C\|\phi\|_{1}\  \delta_D^{\al/2}(x),\quad x\in D.
\end{equation}
\end{lem}
\proof
Starting from \eqref{eq:Gein}, for an arbitrary integer $n\ge1$ we obtain
$$
\phi(x)=\lambda^n\int_D\tilde G^{(n)}_D(x,y)\phi(y)dy,
$$
where $\tilde G^{(1)}_D(x,y):=\tilde G_D(x,y)$ and
$$
\tilde G^{(n+1)}_D(x,y):=\int_D\tilde G_D(x,z)\tilde G^{(n)}_D(z,y)dz.
$$
To estimate $\tG^{(n)}_D$, we  use basic properties of the Bessel potentials, which can be found in \cite[Ch.II.\S4]{MR0143935}.
Recall that for $\al>0$ the  Bessel potential kernel ${\cal G}_{\al}$ is the 
unique, extended-continuous, probability density function on $\Rd$, whose Fourier transform is
 $$
\hat{{\cal G}}_{\al}(\xi)=(1+|\xi|^2)^{-\al/2},\quad \xi\in\bbR^d.
 $$ 
Thus, 
 ${\cal G}_{\al}*{\cal G}_{\beta}={\cal G}_{\al+\beta}$ for all $\al,\beta>0$, see (4.7) of ibid.
If $\al<d$, 
then by 
\cite[(4.2)]{MR0143935}, ${\cal G}_\alpha(x)$ is locally comparable with $|x|^{\al-d}$.
By \eqref{eq:ubRk}  there is a constant $c>0$ such that
$$
\tG_D(x,y)\le c {\cal G}_\alpha(y-x),
 \quad x,y\in \Rd,
 $$
Hence, $\tG^{(n)}(x,y)\le c^n {\cal G}_{n\alpha}(y-x)$,
 which is bounded
if 
$n\alpha>d$, see \cite[(4.2)]{MR0143935} again.
Considering 
such $n$ 
we conclude that $\phi$ is bounded.
The boundary decay  of $\phi$ follows  from \eqref{eq:Gein}, \eqref{eq:egf} and \eqref{oco}. The continuity of $\phi$ is a consequence of the continuity of $\tG_D(x,y)$ for $y\neq x$, and the uniform integrability of the kernel, which stems from \eqref{eq:ubRk}.
\qed

\section{Proof of Theorem~\ref{main-thm}}
\label{sec:pt11}

We say that
$w$ is 
a {\em first integral} of $b$ if 
\begin{equation}\label{eq:dfi} 
\int_D w  b\cdot \nabla \psi = 0,\quad 
\psi\in C^\infty_c(D),
\end{equation}
cf. \eqref{int-first}.
We 
write $w\in{\cal I}_0^{\al}$ if $w \in H^{\al/2}_0(D)$ and $w$ is not equal to $0$ a.e.

Recall that for $A\in\bbR$, the operator $L_A = \Delta^{\alpha/2} + A b\cdot\nabla$
is considered with the Dirichlet exterior condition on $D$, 
i.e.  
it acts on $G_D(L^2(D))$, see Theorem~\ref{thm011307}.
The 
Green operator and  Krein-Rutman eigen-pair of $L_A$ shall be denoted by $\tG_A$ and  $(\lambda_A,\phi_A)$, respectively.
We also recall that $\phi_A\in {\mathcal D}(L_A)$ and
\begin{equation}
\label{022010}
 {\mathcal D}(L_A)={\mathcal D}(\Delta^{\alpha/2})\subset H^1_0(D)\subset H_0^{\alpha/2}(D).
\end{equation}
The proof of  \eqref{081007}  shall be obtained by demonstration  of   lower and upper bounds  
for 
$\lambda_A$ (as $A\to \infty$).

\subsection{Proof of the lower bound (\ref{081007a})}\label{sec:lb}

\begin{prop}
\label{prop010906}
If $f\in {\cal D}(L_A)$, then 
$((-L_A)f,f)={\cal E}^{\alpha}(f,f)$.
\end{prop}
{\em Proof.}
Let $f\in {\cal D}(L_A)$. According to Proposition \ref{prop011009}, $f$ belongs to 
${\mathcal D}(\Delta^{\alpha/2})$ and $ H^1_0(D)$. In addition, $L_Af=\Delta^{\al/2} f+Ab\cdot\nabla f$.
Taking the scalar product of both sides of the  equality against $f$  and using \eqref{antis}
we get
the result, because the second term vanishes.
\qed

According to Proposition \ref{prop010906}, 
\begin{equation}\label{eq:oGp}
\lambda_A=((-L_A)\phi_A,\phi_A)={\cal E}^{\alpha}(\phi_{A},\phi_{A})=(G_D^{-1}\phi_{A},\phi_{A})=\|G_D^{-1/2}\phi_{A}\|_2^2.
\end{equation}
Suppose  
$A_n\to+\infty$, as $n\to \infty$, but $
\lambda_{A_n}$ stay
bounded. By 
\eqref{eq:oGp} and Corollary \ref{cor011206}, 
the sequence $\phi_{A_n}=G_D^{1/2}G_D^{-1/2}\phi_{A_n}$
is pre-compact in $L^2(D)$. 
Suppose that $w$ is 
a weak
limit of $\phi_{A_n}$ 
in $H^{\al/2}_0(D)$, thus a strong limit in $L^2(D)$. We have $\|w\|_2=1$, and for 
$\psi\in C_c^\infty(D)$,
$$
\lambda_{A_n}(\phi_{A_n},\psi)={\cal E}^{\alpha}(\phi_{A_n},\psi)-A_n\int_{D}\phi_{A_n}b\cdot\nabla \psi dx.
$$
Dividing both sides by $A_n$ and passing to the limit, 
we obtain that
$$
\int_{D}w b\cdot\nabla \psi dx=0,
$$
thus $w\in {\cal I}^{\al}$. 
Fatou's lemma and \eqref{eq:oGp} yield
$\liminf\limits_{n\to+\infty}\lambda_{A_n} \ge 
{\cal E}^\alpha(w,w)
$, 
therefore \eqref{081007a} follows.

\subsection{The proof of the upper bound \eqref{091007b}}\label{sec:ub}

The proof of \eqref{091007b} uses "conditioning" of truncations of $w^2$ 
by
the principal eigenfunction inspired by \cite{MR2140256} and \cite{MR2663757}
(see \eqref{071007a} below). Here $w$ is a
first integral in $H^{\alpha/2}_0(D)$. An important part of
the procedure is to prove that  the truncation of $w^2$ is also a first
integral in $H^{\alpha/2}_0(D)$. 
This is true
if $w\in H^1_0(D)$:
suppose that $f:\RR\to\RR$ and $f'$ are bounded and $f(0)=0$.
Let $\diver b=0$ and 
that $\diver(wb)=b\cdot \nabla w=0$ a.e.
Then, a.e. we have
\begin{equation}
\label{020611}
\diver(f(w)b)=f(w)\diver b+f'(w)b\cdot \nabla w
=0\quad\mbox{  on }D,
\end{equation}
thus $f(w)$ is also a first integral of $b$.
However, 
the a.e. differentiability of $w\in H^{\alpha/2}_0(D)$ 
is not  guaranteed for $\alpha<2$. 
In fact, for $w\in H^{\alpha/2}_0(D)$,
the condition $\diver(wb)=0$
is understood 
in the
sense of distributions theory, cf \eqref{int-first}, and the calculations 
in
\eqref{020611} 
may only serve  as a motivation. 
To build up 
tools
for a rigorous proof of the distributional version of \eqref{020611}, 
in Section
\ref{sec5.2.1} 
for a Sobolev-regular, divergence-free  vector field $b:\Rd\to \Rd$,
we consider the incompressible flow  $\{X(t,\cdot),\,t\in\bbR\}$ of 
mappings on $\bbR^d$ 
constructed by DiPerna and Lions in \cite[Theorem III.1]{MR1022305} (the integral curves of $b$). 
The flow is used in Lemma \ref{lem012707} to characterize the first
integrals in $H_0^{\al}(D)$ as 
those
locally
invariant under the flow. 
Then, we are able 
to conclude that the composition of a first integral in $H_0^{\al}(D)$ with a Lipschitz
function is also a first integral in $H_0^{\al}(D)$, see Corollary \ref{prop011007},
and, finally, use the conditioning of $w^2$.

\subsubsection{Flows corresponding to Sobolev regular drifts}

\label{sec5.2.1}
Unless stated otherwise, in this section
we 
consider general $b:\Rd\to \Rd$ such that
$b\in W^{1,1}_{{\rm  loc}}(\Rd)$
and $\diver b=0$  a.e. on $\bbR^d$. 
According to \cite[Theorem III.1]{MR1022305}, there exists a unique 
{a.e. defined jointly Borelian family of mappings
$X(\cdot,x):\RR\to \Rd$ (flow generated by $b$)
with the following properties:
first, for a.e. $x\in \Rd$, the function
$\RR\ni t\mapsto b(X(t,x))$ is continuous and 
\begin{equation}\label{012709a}
X(t,x)=x+\int_0^t b(X(s,x))ds,\quad t\in \RR,
\end{equation}
so that, in particular,
\begin{equation}\label{012709b}
X(0,x)=x\quad \mbox{ and } \quad X(t,X(s,x))=X(t+s,x), \quad s,t\in \RR, 
\end{equation}
second, for all $t\in\bbR$ and Borel measurable 
sets $A\subset \bbR^d$, 
\begin{equation}\label{eq:pLm}
m_d(X(t,A))=m_d(A),
\end{equation}
where $m_d$ is the $d$-dimensional Lebesgue measure, and, third, 
for all $s\in \Rd$,
\begin{equation}\label{012709c}
\lim_{t\to s}\|X(t)-X(s)\|_{L^1(B(0,R))}=0,\quad 
R>0.
\end{equation}
By \eqref{012709b} and \eqref{eq:pLm}, if $f, g$ are 
nonnegative, then 
\begin{equation}\label{eq:switch}
\int_\Rd f(X(t,x))g(x)dx=\int_\Rd f(x)g(X(-t,x))dx.
\end{equation}
We let $p\in[1,+\infty]$, $u_0\in L^p(\bbR^d)$, and define
\begin{equation}\label{eq:defuf}
u(t,x)=u_0(X(t,x)),\quad t\in \RR,\, x\in \bbR^d.
\end{equation}
Note that, $(t,u_0(\cdot))\mapsto u(t,\cdot)$, $t\in \bbR$, defines a group of isometries on $L^p(\bbR^d)$.
If
$v \in C_c^{\infty}(\bbR^{d+1})$, then
\begin{equation}\label{eq:calculus}
v(t,X(t,x))=v(0,x)+\int_0^t \partial_s v(s, X(s,x))ds+
\int_0^t \diver(vb)(X(s,x))ds,
\end{equation}
for a.e. $x\in \Rd$ and all $t\in\bbR$,
because $\diver b=0$ implies 
$b\cdot \nabla v=\diver(v b)$.
Then for all $t\in\bbR$, 
\begin{align*}
&\int_\Rd u(t,x)v(t,x)dx=\int_\Rd u_0(x)v(t,X(-t,x))dx
=\int_\Rd u_0(x)v(0,x)dx\\
&+\int_0^t\int_\Rd u(s,x)\partial_sv(s,x)dxds
-\int_0^t\int_\Rd u(s,x)\diver [v(s,x)b(x)] dxds.
\end{align*}
Let $g\in C_c^\infty(B(0,1))$, 
$g\geq 0$ and $\int_{\bbR^d}g\, dx=1$ (a mollifier).
Define
\begin{equation}\label{eq:mol}
u^{(\eps)}(t,x):=\frac{1}{\eps^d}\int_{\bbR^d}u(t,x-y)g\left(\frac{y}{\eps}\right)dy.
\end{equation}
The function 
is an approximate solution of the transport equation
\begin{equation}
\frac{\partial u}{\partial t}-b\cdot \nabla u=0,\label{transport:eqt}
\end{equation}
 in the following sense: if
\begin{equation}
\label{012807}
r^{(\eps)}(t,x):=\partial_t u^{(\eps)}(t,x)-b(x)\cdot \nabla_xu^{(\eps)}(t,x),
\end{equation}
$q\in[1,+\infty]$, $1/q+1/p=1$ and $b\in W^{1,q}_{{\rm loc}}(\bbR^d)$, then
for all $t\in\bbR$ and finite $R>0$,
\begin{equation}
\label{012807a}
\lim_{\eps\to0+}\int_0^t\|r^{(\eps)}(s,\cdot)\|_{L^q(B(0,R))}ds=0,
\end{equation}
This follows from \cite[Theorem II.1]{MR1022305}.
We shall use $u^{(\eps)}$ to 
characterize the first integrals 
$b\in H^{\alpha/2}_0(D)$
as 
those 
elements of $H^{\alpha/2}_0(D)$,
which are constant along the flow $X(t,\cdot)$.
A word of explanation may be helpful. Suppose that 
 $b\in W^{1,q}_{\rm loc}(\Rd)$, $\diver b=0$ a.e., and $w$ is an $L^p$-integrable
first integral, i.e. 
 $b \cdot\nabla w =0$, as distributions on space-time $\RR\times \Rd$, see \cite[(13)]{MR1022305}.
By \cite[Corollary~II.1]{MR1022305}, there is a unique solution to the transport equation
\eqref{transport:eqt}
with the initial condition $u(0,\cdot)=w$.
The equation is understood in the sense of distributions on space-time, too.
By \cite[Theorem~III.1]{MR1022305}, the solution has the form $(t,x)\mapsto w(X(t,x))$.
However, since $w$ is a first integral, the mapping $(t,x)\mapsto w(x)$ defines another solution.
Thus, by uniqueness, $w(X(t,x))=w(x)$ a.e.  In our case this argument needs to be
slightly modified since the first integral is defined only
on 
$D$ and not on the entire $\Rd$.
In particular, the identity $w(X(t,x))=w(x)$ is bound to hold only for small times $t$. 
\begin{lem}
\label{lem012707}
Let $q=2d/(d+\al)$, $b\in L^\infty(\Rd) \cap W^{1,q}_{{\rm loc}}(\bbR^d)$, and $w\in H^{\al/2}_0(D)$.
Then $w\in {\cal I}^{\al}$ if and only 
for every $\rho>0$ there is
$\kappa>0$ such that 
\begin{equation}
\label{012707}
w(X(t,x))=w(x) \quad\mbox{  for $|t|<\kappa$ and a.e. $x\in D$ satisfying 
$\delta_{D}(x)>\rho$.}
\end{equation}
\end{lem}
\proof  
Suppose that $w\in H_0^{\al/2}(D)$ satisfies
  \eqref{012707} with $\kappa>0$ and let $u_0\in C_c^\infty(D)$. Choose $\rho>0$ 
so that $\delta_{D^c}(x)>\rho$
for 
all $x$ 
in the support
  of $u_0$.
Using \eqref{012707} and 
\eqref{eq:switch},
for all $|t|<\kappa$ we 
obtain,
\begin{equation}
\int_{D}  w(x)u_0(x)dx=\int_{D}  w(X(-t,x))u_0(x)dx=\int_{D}  w(x)u_0(X(t,x))dx.
\label{eq:10711}
\end{equation}
Applying \eqref{eq:calculus}, we rewrite the 
rightmost
side of
\eqref{eq:10711} to obtain,
\begin{equation*}
\int_{D}  w(x)u_0(x)dx=\int_{D}  w(x)u_0(x)dx+ \int_{D}  w(x)\int_0^t b(X(s,x)) \cdot \nabla u_0(X(s,x))ds\,dx.\label{eq:20711}
\end{equation*}
As a result,
$$
0=\int_{D}  w(x)\left\{\frac{1}{t}\int_0^t b(X(s,x)) \cdot \nabla u_0(X(s,x))ds\,\right\}dx.
$$
Letting $t\to0$ we 
see that
$$
0=\int_{D}  w(x)b(x) \cdot \nabla u_0(x)dx.
$$
The limiting passage 
is justified 
by boundedness of $b$, integrability of $w$ (cf. the discussion preceding \eqref{012709a})
and dominated convergence theorem. Since $u_0\in C_c^\infty(D)$ is arbitrary, we 
conclude that $w\in{\cal I}_0^{\al}$.

Conversely, let us assume \eqref{eq:dfi} for some $w\in H^{\al/2}_0(D)$. By Sobolev embedding theorem
\cite[Theorem V.1, p. 119]{MR0290095}, we have that $w\in L^p(D)$,
where $p=2d/(d-\al)$. Let $q=2d/(d+\al)$, $\tilde w=w$ on $D$ and $\tilde w=0$ on $D^c$.
Note that $1/p+1/q=1$. 
Let $\eps>0$ and $t$ be such that 
$$\varepsilon+|t|\|b\|_\infty<{\rm dist}({\rm supp}\, u_0, D^c).$$}
 By \eqref{012709a} and \eqref{eq:mol},
both
$u(t,\cdot)$ and $u^{(\eps)}(t,\cdot)$ are
supported in $D$. 
Using \eqref{012807} and \eqref{eq:dfi} we have
\begin{align*}
\frac{d}{dt}\int_{\bbR^d} \tilde w(x)u^{(\eps)}(t,x)dx
&=\int_{D} w(x)r^{(\eps)}(t,x)dx
+
\int_{D} w(x)b(x)\cdot\nabla_x u^{(\eps)}(t,x)dx\\
&=\int_{D} w(x)r^{(\eps)}(t,x)dx.
\end{align*}
Therefore,
\begin{equation}
\int_{D}  w(x)u^{(\eps)}(t,x)dx-\int_{D} w(x)u_0(x)dx
=\int_0^tds\int_{D} w(x)r^{(\eps)}(s,x)dx.
\label{eq:30711}
\end{equation}
Since $b\in W^{1,q}_{\rm loc}(\bbR^d)$, the remainder $r_\eps(t,x)$ satisfies
\eqref{012807a}. By H\"older inequality,
the right hand side of \eqref{eq:30711} tends to $0$, as $\eps\to0$.
This proves \eqref{012707}
for $|t|<\kappa$.
\qed

As an immediate consequence of Lemma \ref{lem012707} we 
obtain the following.
\begin{cor}
\label{prop011007}
If 
$f$ is Lipschitz on $\bbR$, $f(0)=0$
and 
$w\in{\cal I}_0^\alpha$, then either $f(w)=0$ or
$
f(w)\in {\cal I}^\alpha_0$.
\end{cor}
\proof
By \cite[Theorem 1.4.2 (v)]{MR2778606}, $f(w)\in H^{\alpha/2}_0(D)$.
We apply $f$ to 
\eqref{012707}, and use Lemma \ref{lem012707}.
\qed

\noindent
In particular, we may consider truncations 
of $w$ at the level $N>0$,
$$w_N:=(w\wedge N)\vee (-N).$$
\begin{cor}
\label{cor-x1}
If $N>0$ and $w\in {\cal I}^\alpha_0$, then $w_N\in {\cal I}^\alpha_0$.
\end{cor}

\subsubsection{The upper bound when the drift is defined on entire $\bbR^d$}

\label{sec5.2.2}
We shall first prove the upper bound \eqref{091007b} under
the assumptions 
that $b$ is bounded and of zero divergence on the whole of $\Rd$, 
and $b\in W^{1,2d/(d+\alpha)}_{{\rm loc}}(\bbR^d)$.
\begin{prop}
\label{prop021007}
Suppose that $A\in \bbR$, $\eps>0$,
$w\in{\cal I}^{\al}_0$
and $w$ is bounded.
Then,
\begin{eqnarray}
\label{061007}
&&
\lambda_A\int_D\frac{\phi_A(z) w^2(z)}{\phi_A(z)+\eps}dz\\
&&
=\frac{1}{2}\int_{D^2}\left\{\frac{w^2(x)}{\phi_A(x)+\eps}-\frac{w^2(y)}{\phi_A(y)+\eps}\right\} 
\frac{\phi_A(x)-\phi_A(y)}{|x-y|^{d+\al}}\nu(x-y)dxdy.\nonumber
\end{eqnarray}
\end{prop}
\proof 
We denote 
$\psi=w^2/(\phi_A+\eps)$.
By \cite[Theorem 1.4.2 (ii), (iv)]{MR2778606}, $\psi\in H_0^{\alpha/2}(D)$. Also, $\log(\phi_A+\eps)-\log \eps\in  H_0^{1}(D)$.
Considering 
\eqref{011109b}, we observe
that the right hand side of \eqref{061007} 
equals ${\cal E}_{\al}(\phi_A,\psi)$.
By Proposition \ref{prop011009}, the left hand side of \eqref{061007} is
$$-\int_DL_A\phi_A(z)\psi(z)dz
=-\int_D\Delta^{\al/2}\phi_A(z)\psi(z)dz-\int_Db(z)\cdot \nabla\phi_A(z)\psi(z)dz.
$$
However, the second term on the right hand side vanishes, 
because
\begin{equation*}\label{eq:w2cp}
-\int_Dw^2(z)b(z)\cdot \nabla\log[\phi_A(z)+\eps]dz=0,
\end{equation*}
and $w^2=w^2_N$, for a sufficiently large $N$, is a
first integral by virtue of Corollary \ref{cor-x1}. 
Thus \eqref{061007} follows from
\eqref{eq:GfD}.
\qed

The following elementary  identity holds for 
functions 
$u,v$,
\begin{eqnarray}
\label{071007}
&&[u(x)-u(y)]^2+u^2(x)\frac{v(y)-v(x)}{v(x)}+u^2(y)\frac{v(x)-v(y)}{v(y)}\\
&&
=v(x)v(y)\left[\frac{u(x)}{v(x)}-\frac{u(y)}{v(y)}\right]^2\nonumber.
\end{eqnarray}
Let 
$N,\eps>0$,
$u=w_N$, and $v=\phi_A+\eps$.
By \eqref{071007},
\begin{equation}
\label{071007a}
[w_N(x)-w_N(y)]^2\ge \left\{\frac{w^2_N(x)}{\phi_A(x)+\eps}-\frac{w^2_N(y)}{\phi_A(y)+\eps}\right\}\left[\phi_A(x)-\phi_A(y)\right].
\end{equation}
Multiplying both sides 
by $\nu(x-y)dxdy$ and integrating over $D^2$,
we obtain
$$
2{\cal E}^{\alpha}(w_N,w_N)\ge \int_{D^2}\left\{\frac{w^2_N(x)}{\phi_A(x)+\eps}-\frac{w^2_N(y)}{\phi_A(y)+\eps}\right\}\left[\phi_A(x)-\phi_A(y)\right]\nu(x-y)dxdy.
$$
Using 
\eqref{061007}, we 
see that
$$
{\cal E}^{\alpha}(w_N,w_N)\ge \lambda_A \int_{\bbR^d}\frac{w^2_N(x)\phi_A(x)dx}{\phi_A(x)+\eps}.
$$
Letting $\eps\to0+$, we conclude that
${\cal E}^{\alpha}(w_N,w_N)\ge\lambda_A \|w_N\|_2^2$.
Letting $N\to+\infty$  
and using \cite[part (iii) of Theorem 1.4.2]{MR2778606}, 
we obtain
${\cal E}^{\alpha}(w,w)\ge\lambda_A\|w\|_2^2$.
This proves \eqref{091007b}.

\subsubsection{
The upper bound when the drift is defined only on $D$
}\label{secmain-thm} 

\label{sec5.2.3}

\label{sec6}

Suppose that $b\in W^{1,q}(D)$, with $q=2d/(d+\al)$, and $\diver b(x)\equiv0$ 
on $D$. Here, as usual, $D$
is a bounded domain with the $C^{1,1}$ class boundary. By the discussion in Section~\ref{sec:gener},
$\bbR^d\setminus\bar D$ has finitely many, say, $N+1$ connected
 components. Denote them by ${\cal O}_0,\ldots, {\cal O}_N$, 
and assume that $\infty$ belongs to the compactification of ${\cal O}_0$.
We start with the following extension result. 
\begin{prop}
\label{prop013110}
There exist $\delta>0$ and $\tilde b\in W^{1,q}(\bbR^d)$ such that 
\begin{equation}
\label{023110}
\diver \tilde b\in L^\infty(\bbR^d),\qquad 
\tilde b
=b \quad \mbox{on $D$,}
\end{equation}
$\tilde b$ has compact support and is smooth outside of 
$D_{\delta}:=\set{x\in \Rd:
{\rm dist}(x,D)<\delta}$.
\end{prop}
\proof
According to the results of Section 4 of   \cite{MR1809290} we can
find $h_{ij}\in W^{2,q}(\bbR^d)$, supported in $D_{\delta}$ and such that
$h_{ij}(\cdot)=-h_{ji}(\cdot)$ for $i,j=1,\ldots,d$, numbers $\la_1,\ldots,\la_N\in\bbR$, and points $y^{(1)}\in{\cal O}_1,\ldots,y^{(N)}\in{\cal O}_N$ such that
 $$
b_i(x)=\sum_{j=1}^dh_{ji}(x)+\sum_{j=1}^N\la_j\frac{x_i-y_i^{(j)}}{|x-y^{(j)}|^d},\quad \mbox{for } x\in D,\,i=1,\ldots,d.
$$
Here $b=(b_1,\ldots,b_d)$. We denote by $\hat b_i$ the right hand side of the above equality. The
field $\hat b$
is incompressible on $\Rd\setminus \set{y^{(1)},\ldots,y^{(d)}}$.
We choose $\delta>0$ so that ${\rm dist}(y^{(j)},D)>4\delta$
for $j=1,\ldots,N$. 
Let $\phi
\in C^\infty$ 
equal to $1$ on $D_{\delta/2}$ and 
vanish outside of $D_\delta$.
The field $\tilde b:=\hat b \phi$ has the desired properties.
\qed

We consider the
flow
$
X
$ of measurable mappings 
generated by $\tilde b$, which satisfies \eqref{012709a},
\eqref{012709b} and \eqref{012709c}, see Theorem III.2 of \cite{MR1022305}. Condition \eqref{eq:pLm} 
needs to be
modified 
as follows: there exists $C>0$, such that 
\begin{equation}
\label{073110}
e^{t/C}m_d(A)\le 
m_d(X(t,A))\le e^{Ct}m_d(A), \quad A\subset \Rd,\; t\in \RR.
\end{equation}
Let $g
$ be a mollifier, let
$$
b^{(\eps)}(x):=\frac{1}{\eps^d}\int_{\bbR^d}\tilde b(x-y)g\left(\frac{y}{\eps}\right)dy,
$$
and let 
$
X_\eps
$ be the flow generated by  $b^{(\eps)}
$.
It has been shown in Section 3 of \cite{MR1022305} that
\begin{equation}
\lim_{\eps\to0}\sup_{t\in[-T,T]}\int_{B(0,N)}|\phi(X(t,x))-\phi(X_\eps(t,x))|\wedge 2^Ndx=0\label{eq:2.31.10}
\end{equation}
for all $T,N>0$ and measurable functions $\phi:\bbR^d\to\bbR$.
We shall need the
following 
modification of Lemma \ref{lem012707}.
\begin{lem}
\label{lem012707a}
If
$b\in L^\infty(D)\cap W^{1,2d/(d+\al)}(D)$ then, 
the conclusion of Lemma \ref{lem012707} holds
for the flow
$
X
$ 
generated by  $\tilde b$.
\end{lem}
\proof Suppose that $w\in {\cal I}^{\al}_0$. 
Since
$\tilde b=b$ on $D$,
we can repeat the proof of the respective part of Lemma \ref{lem012707}.
Namely, keeping the notation from that lemma, 
for every $u_0\in C_c^\infty(D)$ 
we have $\ka>0$ such that
\begin{equation}
\label{022807a}
\int_{D}  w(x)u_0(X(-t,x))dx=\int_{D} w(x)u_0(x)dx,\quad
|t|<\kappa.
\end{equation}
To complete the proof of \eqref{012707}
it suffices 
to conclude that $\ka>0$ can be so adjusted 
that
\begin{equation}
\int_{D}  w(x)u_0(X(-t,x))dx
=\int_{D} u_0(x)w(X(t,x))dx
,\quad
|t|<\kappa.
\label{eq:1.31.10}
\end{equation}
This part cannot be guaranteed directly from the definition of the
flow, as the extended field $\tilde b$ needs not be 
divergence-free.
Equality \eqref{eq:1.31.10} holds however, when $X(t)$ is replaced by
$X_\eps(t)$ for a sufficiently small $\eps>0$. Indeed, by the Liouville
theorem, the Jacobian
${\cal J} X_\eps(t,x)$  of  $ X_\eps(t,x)$ satisfies 
\begin{eqnarray*}
&&
\frac{d}{dt}{\cal J} X_\eps(t,x)=\diver b^{(\eps)}(
X_\eps(t,x)) {\cal J} X_\eps(t,x),\\
&&
 {\cal J} X_\eps(0,x)=1,\quad
t\in \RR,\, x\in\Rd.
\end{eqnarray*}
Since $\diver b^{(\epsilon)}=0$ in an open neighborhood of $\bar D$
we conclude that ${\cal J} X_\eps(t,x)\equiv 1$ on $D$ for 
(sufficiently small) 
$|t|<\ka$.
Since
$w$ and $u_0$ are supported in $D$,
\begin{equation}
\int_{D} u_0(x)w(X_\eps(t,x))dx=
\int_{D}  w(x)u_0(X_\eps(-t,x))dx,\quad
|t|<\kappa.
\label{eq:1.31.10eps}
\end{equation}
 Letting $\eps\to0$ and using \eqref{eq:2.31.10} we 
obtain \eqref{eq:1.31.10}. The rest of the proof
follows 
that of Lemma \ref{lem012707}.
\qed

Having established Lemma \ref{lem012707a} we proceed with the proofs of Corollaries \ref{prop011007} and \ref{cor-x1} and Proposition \ref{prop021007} with no alterations. These results 
yield \eqref{091007b}.

\begin{example}
\rm
We consider the principal eigen-pair, say $(\lambda_0,\phi_0)$, of $\Delta^{\alpha/2}$ for the unit ball in $\bbR^d$. By rotation invariance of $\Delta^{\alpha/2}$ and uniqueness, $\phi_0$ is  a smooth radial function in the ball. 
For $i,j=1,\ldots,d$ we take  radially symmetric functions 
$h_{ij}(|x|^2)\in   C_c^\infty(\bbR^d)$, such that $h_{ij}=-h_{ji}$.
Let $b_i=\sum_{j=1}^d\partial_jh_{ji}$, $i=1,\ldots,d$.
The 
vector field $b=(b_1,\ldots,b_d)$ is of  zero
divergence and 
tangent to the spheres $|x|=r$ for all $r>0$. Indeed, since $h'_{ij}=-h'_{ji}$,
$$
b(x)\cdot x=\sum_{i,j=1}^d\partial_{j}h_{ji}(|x|^2)x_i=2\sum_{i,j=1}^dh_{ji}'(|x|^2)x_ix_j\equiv0.
$$
As a result
$b(x)\cdot \nabla w(x)= 0$ for any $C^1$ smooth radially symmetric function $w(\cdot)$.
Thus, we conclude that $\phi_0$ is the principal eigenfunction of $L_A$ for every $A$. We have $\lambda_A\equiv \lambda_0$, and ${\cal E}^\alpha(w,w)$ attains its infimum, $\lambda_0$, at $w=\phi_0$, see Section~\ref{sec:lb}. 
In passing we note that the considered limiting eigenproblems are essentially different for different values of $\alpha$, in accordance
with the fact that the "escape rate"
$x\mapsto \int_{D^c} \nu(y-x)dy$ of the isotropic $\alpha$-stable L\'evy processes from $D$ depends on $\alpha$. We refer the interested reader to \cite{MR2974318}, \cite{MR2876409}, \cite{MR2217951} for more information on the eigenproblem of $\Delta^{\alpha/2}$, see also  \cite{MR2158176}.
\end{example}

\subsection{Existence of a minimizer}

\label{proof-coro} 
To complete the proof of Theorem~\ref{main-thm}
we only need to explain the attainability of infimum
appearing 
on the right hand side of \eqref{081007}. 
This is done in the following.
\begin{lem}
\label{lem1}
If $\set{w\in {\cal I}_0^{\alpha/2}:\,\|w\|_2=1}\neq \emptyset$, then
$w\mapsto {\cal E}^\alpha(w,w)$ attains 
its infimum 
on the set.
\end{lem}
\proof
If $e_*<\infty$ is the infimum,
then we can choose functions $w_n$ in the set given in the statement of the lemma, such that ${\cal E}^\alpha(w_n,w_n)\to e_*$ 
and, by choosing a subsequence, that $w_n$ weakly converge to $w_*\in H_0^{\alpha/2}(D)$. 
Since $\set{w_n}$ is precompact in $L^2(D)$, we may further assume that
$w_n\to w$ in $L^2(D)$ and a.e. 
This implies that $\|w_*\|=1$, $w_*$ is a first integral and,
by \eqref{011109b} and Fatou's lemma, that ${\cal E}^\alpha(w_*,w_*)\le e_*$. From the definition of $e_*$, we conclude that 
equality actually occurs.
\qed

\section*{Acknowledgements}
Krzysztof Bogdan thanks the Department of Mathematics at Stanford University for hospitality during his work on the paper, Lenya Ryzhik for discussions on incompressible flows and Michael Frazier for a discussion on Schur's test.

\def\cprime{$'$}

\end{document}